\documentclass[11pt]{article}
\usepackage{ graphics}

\newcommand{\mathsym}[1]{{}}
\newcommand{\unicode}[1]{{}}
\usepackage{latexsym,xy,eucal,mathrsfs,graphs}
\usepackage{setspace}
\usepackage{amsthm}
\usepackage{pstricks}
\usepackage{amssymb}
\usepackage{amsmath}
\usepackage{indentfirst}
\usepackage{epsfig}
\usepackage{indentfirst}
\usepackage{amscd}
\usepackage{multirow}

\newtheorem{theorem}{Theorem}[section]
\newtheorem{proposition}[theorem]{Proposition}
\newtheorem{lemma}[theorem]{Lemma}

\newtheorem{corollary}[theorem]{Corollary}

\newtheorem{result}[theorem]{Result}
\theoremstyle{definition}
\newtheorem{defn}[theorem]{Definition}
\newtheorem{example}[theorem]{Example}

\usepackage{simplemargins}
\setleftmargin{0.88in}
\setrightmargin{0.88in}
\settopmargin{0.817in}
\setbottommargin{0.83in}

 \font\conj =
msbm10 at 12pt  

\def\text#1{\mathop{{\rm{#1}}}\nolimits}
\def\mr{\mathop{{\rm{mr}}}\nolimits}
\def\msr{\mathop{{\rm{msr}}}\nolimits}

\def\cc{\mathop{{\rm{cc}}}\nolimits}

\def\d{\mathop{{\rm{d}}}\nolimits}
\def\sp{\mathop{{\rm{Span}}}\nolimits}

\def\Re{{\mbox{\conj R}}}
\def\N{{\mbox{\conj N}}}

\def\F{{\mbox{\conj F}}}
\def\C{{\mbox{\conj C}}}

\def\dis{\displaystyle}
\def\1{\'{\i}}

\title{ A Proof   for   Delta Conjecture}

\author{Pedro D\'{\i}az Navarro\thanks{Escuela de Matem\'atica, Universidad de Costa Rica}}
\date{Junio , 2018}
\begin{document}
\maketitle

\begin{abstract}
\noindent
 By finding orthogonal representation for a family of   simple  connected  called $\delta$-graphs it is possible to show that $\delta$-graphs satisfy delta conjecture. An extension of the argument to graphs of the form $\overline{P_{\Delta(G)+2}\sqcup G}$ where $P_{\Delta(G)+2}$ is a path and $G$ is a simple connected graph it is possible to find an orthogonal representation of $\overline{P_{\Delta(G)+2}\sqcup G}$ in $\mathbb{R}^{\Delta(G)+1}$. As a consequence we prove delta conjecture.\end{abstract}

\\
\noindent{{\bf DOI:} 05C50,05C76 ,05C85 ,68R05 ,65F99,97K30.}

\section{Introduction}

A {\it graph} $G$ consists of a set  of vertices $V(G)=\{1,2,\dots,n\}$ and a set of edges $E(G)$, where an edge is defined to be an unordered  pair of vertices. The {\it order} of $G$, denoted $\vert G\vert $,  is the cardinality  of $V(G)$. A graph is {\it simple} if it has no multiple  edges  or loops. The {\it complement } of a graph $G(V,E)$ is the graph $\overline{G}=(V,\overline{E})$, where $\overline{E}$ consists of all those edges of the complete  graph $K_{\vert G\vert}$ that are not in $E$.

 A matrix $A=[a_{ij}]$ is {\it combinatorially symmetric} when $a_{ij}=0$ if and only if $a_{ji}=0$. We say that  $G(A)$  is the graph of a combinatorially symmetric matrix $A=[a_{ij}]$ if $V=\{1,2,\dots,n\}$ and $E=\{\{i,j\}: a_{ij}\ne0\}$ . The main diagonal entries of $A$ play no role in determining $G$. Define $S(G,\F)$ as the set of all $n\times n$ matrices that are real symmetric if $\F=\Re$ or complex Hermitian if $\F=\C$ whose graph is $G$. The sets $S_+(G,\F)$ are the corresponding subsets of positive semidefinite (psd) matrices. The smallest possible rank  of any matrix $A\in S(G,\F)$  is the {\it minimum rank} of $G$, denoted $\mr(G,\F)$, and the smallest possible rank of any matrix $A\in S_+(G,\F)$  is the {\it minimum semidefinite rank} of $G$, denoted $\mr_+(G)$ or $\msr(G)$.

In 1996, the minimum rank among real symmetric matrices with a given graph was studied  by  Nylen \cite{PN}. It gave rise to the area of minimum rank problems which led to the study of minimum rank among complex Hermitian matrices and positive semidefinite matrices associated with a given graph. Many results can be  found for  example  in  \cite{FW2, VH, YL, LM, PN}.

During  the   AIM workshop of 2006 in Palo Alto, CA, it  was  conjectured   that for any   graph $G$  and  infinite field  $F$, $\mr(G,\F)\le |G|-\delta(G)$ where $\delta(G)$  is the minimum  degree  of  $G$.  It  was  shown  that  for if  $\delta(G)\le 3$ or  $\delta(G)\ge |G|-2$    this inequality  holds.  Also  it  can be verified that  if  $|G|\le 6$ then  $\mr(G,F)\le |G|-\delta(G)$.  Also  it  was   proven  that any  bipartite graph  satisfies  this conjecture. This  conjecture   is called  the {\it Delta Conjecture}. If we   restrict  the study  to  consider  matrices  in $S_+(G,\F)$  then  delta conjecture  is  written  as $\msr(G)\le |G|-\delta(G)$. Some  results on delta conjecture  can be found in \cite{AB,RB1, SY1, SY} but  the   general problem  remains  unsolved. In  this paper,  by  using  a  generalization  of  the argument  in  \cite{PD},  we  give  an  argument  which  prove  that  delta  conjecture  is  true   for  any  simple and  connected  graph  which means  that delta  conjecture  is  true.

\section{Graph Theory Preliminaries}
\addtocontents{toc}{\vspace{-15pt}}
In  this section   we give   definitions and results from  graph theory which    will  be used in  the remaining  chapters. Further details  can be found  in \cite{BO,BM, CH}.

A {\bf graph} {$G(V,E)$} is a pair {$(V(G),E(G)),$} where {$V(G)$} is the set of vertices and {$E(G)$} is the set of edges together  with an  {\bf  incidence  function} $\psi(G)$ that associate with  each edge  of  $G$ an  unordered  pair  of (not necessarily  distinct) vertices  of  $G$. The {\bf order} of {$G$}, denoted {$|G|$}, is the number of vertices in {$G.$} A graph is said to be {\bf simple} if it has no loops or multiple edges. The {\bf complement} of a graph {$G(V,E)$} is the graph {$\overline{G}=(V,\overline{E}),$} where {$\overline{E}$} consists of all the edges that are not in {$E$}.
A {\bf  subgraph} {$H=(V(H),E(H))$} of {$G=(V,E)$} is a graph with {$V(H)\subseteq V(G)$} and {$E(H)\subseteq E(G).$} An {\bf induced subgraph} {$H$} of {$G$}, denoted G[V(H)], is a subgraph with {$V(H)\subseteq V(G)$} and {$E(H)=\{\{i,j\} \in E(G):i,j\in V(H)\}$}. Sometimes  we  denote the  edge $\{i,j\}$ as $ij$.

We  say  that  two  vertices of a graph $G$  are {\bf adjacent}, denoted  $v_i\sim v_j$,   if  there is an edge $\{v_i,v_j\}$ in  $G$.  Otherwise  we say  that the  two  vertices $v_i$ and $v_j$ are {\bf non-adjacent}  and  we denote this  by $v_i \not\sim v_j$.  Let {$N(v)$} denote the set of vertices that are adjacent to the vertex {$v$} and let {$N[v]=\{v\}\cup N(v)$}. The {\bf degree} of a vertex {$v$} in {$G,$} denoted {$\d_G(v),$} is the cardinality of {$N(v).$} If {$\d_G(v)=1,$} then {$v$} is said to be a {\bf pendant} vertex of {$G.$} We use {$\delta(G)$} to denote the minimum degree of the vertices in {$G$}, whereas {$\Delta(G)$} will denote the maximum degree of the vertices in {$G$}.

 Two  graphs $G(V,E)$ and $H(V',E')$  are  identical  denoted  $G=H$, if  $V=V',  E=E'$, and $\psi_G=\psi_H$  . Two  graphs $G(V,E)$ and $H(V',E')$ are {\bf isomorphic}, denoted  by  $G\cong H$, if  there exist bijections $\theta:V\to V'$  and $\phi: E\to  E' $ such  that $\psi_G(e)=\{u,v\}$  if  and  only if  $\psi_H(\phi(e))= \{\theta(u), \theta(v)\}$.

 A {\bf complete graph}  is a simple graph in which the vertices are pairwise adjacent.
  We will use {$nG$} to denote {$n$} copies of a graph {$G$}. For example, $3K_1$  denotes three  isolated vertices $K_1$ while {$2K_2$} is the graph given  by   two  disconnected  copies  of $K_2$.

 A {\bf path} is a list of distinct vertices in which successive vertices are connected by edges. A path on {$n$} vertices is denoted by {$P_n.$} A graph {$G$} is said to be {\bf connected} if there is a path between any two vertices of {$G$}. A {\bf cycle} on {$n$} vertices, denoted {$C_n,$} is a path such that the beginning vertex and the end vertex are the same. A {\bf tree} is a connected graph with no cycles. A graph $G(V,E)$ is  said  to be {\bf chordal}  if it has no induced cycles $C_n$ with $n\ge 4$.
 A  {\bf component}  of a graph $G(V,E)$ is  a maximal connected  subgraph. A  {\bf cut vertex}   is  a vertex  whose deletion  increases  the number of  components.

 The {\bf union}  $G\cup G_2$ of  two  graphs $G_1(V_1,E_1)$ and $G_2(V_2,G_2)$  is  the union  of  their  vertex  set  and  edge  set,  that is $G\cup G_2(V_1\cup V_2,E_1\cup E_2$. When  $V_1$ and $V_2$ are disjoint their union  is called  {\bf  disjoint union} and  denoted $G_1\sqcup G_2$.


\section{The Minimum  Semidefinite Rank  of  a Graph}
In  this section  we will establish   some of  the   results for  the minimum  semidefinite rank ($\msr$)of a graph $G$  that  we  will be using in the subsequent chapters.

A {\bf positive  definite} matrix  $A$ is an Hermitian   $n\times n$ matrix such  that $x^\star A x>0$  for all nonzero  $x\in \C^n$. Equivalently,  $A$  is a  $n\times n$ Hermitian positive definite matrix  if and  only  if   all the  eigenvalues of $A$ are positive (\cite{RC}, p.250).

A $n\times n$ Hermitian matrix  $A$ such  that $x^\star A x\ge 0$  for all $x\in \C^n$ is  said  to be  {\bf positive  semidefinite (psd)}. Equi\-va\-lently,   $A$ is a  $n\times n$ Hemitian positive  semidefinite matrix if and  only  if  $A$  has  all   eigenvalues nonnegative (\cite{RC}, p.182).

If $\overrightarrow{V}=\{\overrightarrow{v_1},\overrightarrow{v_2},\dots, \overrightarrow{v_n}\}\subset \Re^m$  is a set of  column vectors  then  the  matrix
$ A^T A$, where $A= \left[\begin{array}{cccc}
  \overrightarrow{v_1} & \overrightarrow{v_2} &\dots& \overrightarrow{v_n}
\end{array}\right]$
and $A^T$  represents  the  transpose matrix of  $A$, is a psd matrix  called  the  {\bf Gram matrix} of $\overrightarrow{V}$. Let $G(V,E)$  be a graph  associated  with  this Gram matrix. Then  $V_G=\{v_1,\dots, v_n\}$ correspond to  the set of  vectors in $\overrightarrow{V}$ and  E(G) correspond to  the nonzero inner products  among  the  vectors  in $\overrightarrow{V}$. In this  case $\overrightarrow{V}$  is  called an  {\bf orthogonal representation} of $G(V,E)$ in $\Re^m$. If  such  an  orthogonal  representation  exists  for  $G$ then $\msr(G)\le m$.

Some    results about  the minimum semidefinite  rank of a graph  are  the  following:

\begin{result}\cite{VH}\label{msrtree}
If  $T $ is a tree  then $\msr(T)= |T|-1$.
\end{result}
\begin{result}\cite{MP3}\label{msrcycle}
 The cycle $C_n$ has minimum semidefinite rank $n-2$.
\end{result}


\begin{result}\label{res2}
 \cite{MP3} \ If a connected  graph $G$  has a pendant  vertex $v$, then $\msr(G)=\msr(G-v)+1$ where $G-v$  is obtained as an induced subgraph of $G$ by  deleting $v$.
\end{result}

\begin{result} \cite{PB} \label{OS2}
 If {$G$} is a connected, chordal graph, then {$\msr(G)=\cc(G).$}
\end{result}

\begin{result}\label{res1}
\cite{MP2}\ If a graph $G(V,E)$ has a cut vertex, so that $G=G_1\cdot G_2$, then  $\msr(G)= \msr(G_1)+\msr(G_2)$.
\end{result}

\section{Delta-Graphs and  the Delta Conjecture}

In  \cite{PD}  is  is  defined   a   family of  graphs  called $\delta$-graphs  and  show  that  they satisfy  the  delta conjecture.

\begin{defn}\label{ccpg}
Suppose  that $G=(V,E)$  with $|G|=n \ge 4$  is simple  and  connected such  that  $\overline{G}=(V,\overline{E})$  is also  simple  and  connected. We  say  that  $G$ is a {\bf $\mathbf{\delta}$-graph} if  we  can  label  the vertices of $G$ in such a way that
\begin{enumerate}
  \item[(1)] the  induced    graph    of  the  vertices  $v_1,v_2,v_3$  in $G$ is  either $3K_1$  or  $K_2 \sqcup K_1$,  and
  \item[(2)] for $m\ge 4$,  the vertex  $v_m$  is  adjacent   to all   the  prior  vertices $v_1,v_2,\dots,v_{m-1}$  except  for at most $\dis{\left\lfloor\frac{m}{2}-1\right\rfloor}$ vertices.
  \end{enumerate}
 \end{defn}
  A second  family  of  graphs also   defined  in  \cite{PD} contains  the  complements  of $\delta$-graphs.
\begin{defn} Suppose   that  a graph $G(V,E)$  with $|G|=n \ge 4$  is  simple and  connected   such  that  $\overline{G}=(V,\overline{E})$  is also   simple and  connected.  We  say  that  $G(V,E)$  is a {\bf C-$\mathbf{\delta}$  graph}  if    $\overline{G}$   is a  $\delta$-graph.

In other  words,  $G$   is a  {\bf C-$\mathbf{\delta}$ graph} if  we can  label   the  vertices   of  $G$  in  such a  way   that
 \begin{enumerate}
 \item[ (1)] the  induced  graph  of  the vertices  $v_1,v_2,v_3$ in  $G$ is  either  $K_3$ or $P_3$,  and
 \item[(2)]  for  $m\ge 4$,  the vertex $v_m$  is adjacent  to at most  $\dis{\left\lfloor\frac{m}{2}-1\right\rfloor}$ of  the  prior  vertices $v_1,v_2,\dots,v_{m-1}$.
\end{enumerate}
\end{defn}
\begin{example}\label{examplecp}
The   cartesian product $K_3\square P_4$  is a C-$\delta$ graph and  its  complement  is  a  $\delta$-graph. By labeling as  the   following picture  we  can  verified  the definition  for  both  graphs.
\newpage
\begin{center}
\includegraphics[height=40mm]{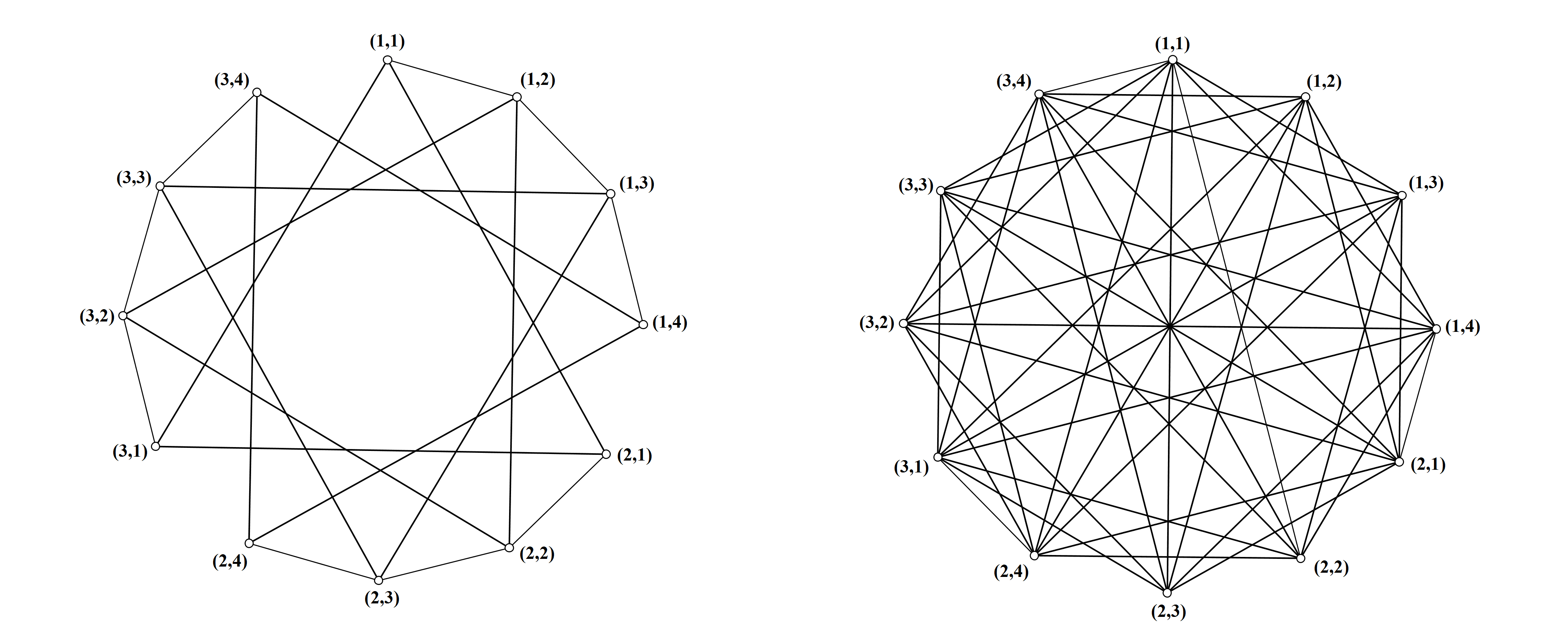}
 \end{center}
\vspace{-0.1in}\begin{figure}[h]
\centering
\caption{The  Graph $K_3\square P_4$  and  its  complement $\overline{K_3\square P_4}$ }\label{fig3.1}
\end{figure}
Note that we  can label  the vertices  of $K_3\square P_4$ clockwise $v_1=(1,1),v_2= (1,2),v_3=(1,3),\dots v_{12}=(3,4)$. The  graph  induced   by  $v_1,v_2,v_3$  is $P_3$. The    vertex $v_4$ is   adjacent to a prior vertex   which  is  $v_3$  in  the induced  subgraph of  $K_3\square P_4$ given  by $\{v_1,v_2,v_3,v_4\}$. Also,   the  vertex $v_5$ is adjacent only  to  vertex $v_1$  in  the  induced  subgraph of  $K_3\square P_4$ given  by $\{v_1,v_2,v_3,v_4, v_5\}$. Continuing   the process trough vertex $v_12$  we conclude   that  $K_3\square P_4$   is a  C-$\delta$  graph. In  the  same  way  we conclude  that  its  complement $\overline{K_3\square P_4}$  is a $\delta$-graph.
\end{example}
\begin{lemma}\label{lem2} Let {$G(V,E)$} be a $\delta$-graph. Then   the induced  graph of {$\{v_1,v_2,v_3\}$} in {$G$} denoted  by $H$ has an  orthogonal  representation  in {$\Re^{\Delta(\overline{G})+1}$} satisfying the  following  conditions:
\begin{enumerate}
\item [(i)] the  vectors   in  the orthogonal  representation  of  $H$ can be chosen  with nonzero coordinates, and
\item [(ii)]\label{L1}$\overrightarrow{v}\not\in \sp(\overrightarrow{u})$ for each pair   of  distinct  vertices  $u,v$ in  $H$.
\end{enumerate}
\end{lemma}
\begin{theorem} \label{main} Let $G(V,E)$ be a $\delta$-graph then
$$
\msr(G)\le\Delta(\overline{G})+1=|G|-\delta(G)\label{mrsineq1}
$$
\end{theorem}
The  proof  of  these  two  results  can  be  found  in  \cite{PD} and \cite{PD1}. The  argument  of  the  proof  is  based on   the construction  of a  orthogonal  representation  of  pairwise linear  independent vectors  for  a $\delta$  graph $G$  at $\Re^{\overline{G})+1}$.  Since  $\msr(G)$  is   the minimum dimension  in which   we can  get  an  orthogonal representation    for a  simple  connected  graphs  the  result  is a  direct  consequence of this  construction.

\section{A survey of  $\delta$-graphs  and upper bounds  their  minimum Semidefinite rank}
The  theorem \ref{main}   give us a  huge family  of  graph  which  satisfies  delta conjecture. Since,  the   complement    of   a  C-$\delta$ graphs  is ussually  a $\delta$-graph,  it  is  enough  to identify  a  C-$\delta$-graph    and  therefore  we know   that   its  complement  is a $\delta$-graph satisfying  delta  conjecture if  it  is  simple and  connected.

Some  examples   of  C-$\delta$ graphs that  we can  find  in  \cite{PD} are  the  Cartesian Product $K_n\square P_m,n\ge  3,  m\ge  4$, Mobi\"us Lader $ML_{2n}, n\ge 3$, Supertriangles $Tn, n\ge 4$,  Coronas $S_n\circ P_m, n\ge  2 , m\ge 1$ where  $S_n$  is a  star  and  $P_m$ a path, Cages like Tutte's (3,8) cage,  Headwood's (3,6) cage   and many others, Blanusa Snarks  of type  $1$ and  $2$  with  $26, 34$, and $42$ v\'ertices, and  Generalized Petersen Graphs $Gp1$ to  $Gp16$.

\subsection{Upper bounds   for  the  Minimum semidefinite rank of  some  families  of  Simple connected graphs}

 From the definition  of  C-$\delta$ graph   and the Theorem \ref{main} we  can obtain  upper  bounds   for  the graph  complement of a C-$\delta$  graphs. It is enough to  label  the  vertices  of $G$  in  such a way  that  the labeled sequence  of  vertices  satisfies  the definition. That is,  if  we  start   with  the induced  graph of  $ \{v_1,v_2,v_3\}$,  the  newly added  vertex $v_m$  is adjacent  to at most $\lfloor\frac{m}{2}-1\rfloor $  of  the prior  vertices $v_1,v_2,\dots,v_{m-1}$.   Then  $G$  is a C-$\delta$  graph  and  its graph  complement $\overline{G}(V,\overline{E})$ will  have an orthogonal representation   in $\Re^{\Delta(G)+1}$ any time it is  simple  and  connected.
 
 In  order   to  show  the  technique  used  in  the proved  result   consider  the following  examples
 
 \begin{example}\label{upper1}

If  $G$ is  the Robertson's (4,5)-cage on 19 vertices  then  it is  a 4-regular C-$\delta$ graph. Since $\Delta(G)=4$, the $\msr(G)\le 5$.  To  see  this is a C-$\delta$ graph it  is enough  to label its  vertices  in the  way  shown  in the next  figure:
\begin{center}
\includegraphics[height=50mm]{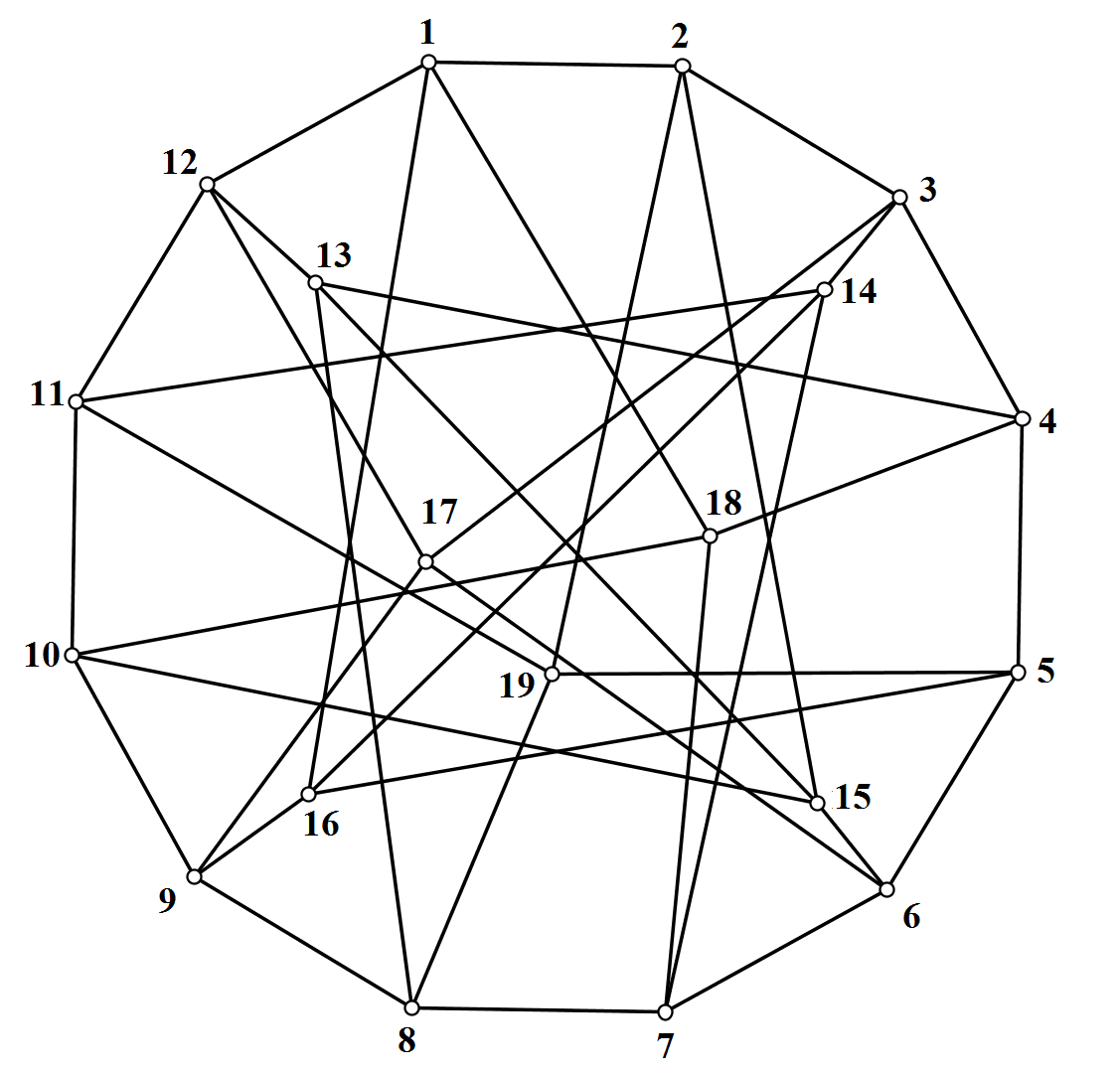}
 \end{center}
\vspace*{-0.1in}\begin{figure}[h]
\centering
{Figure B.2 Robertson's (4,5)-cage (19 vertices)}
\label{figA.1.2}
\end{figure}
\end{example}
\begin{example}\label{upper2}
If  $G$ is  the platonic  graph   Dodecahedron  then  it is  a 3-regular C-$\delta$ graph. Since $\Delta(G)=3$, the $\msr(\overline{G})\le 4$.  To  see  this is a C-$\delta$ graph it  is enough  to label its  vertices  in the  way  shown  in the next  figure:
\begin{center}
\includegraphics[height=60mm]{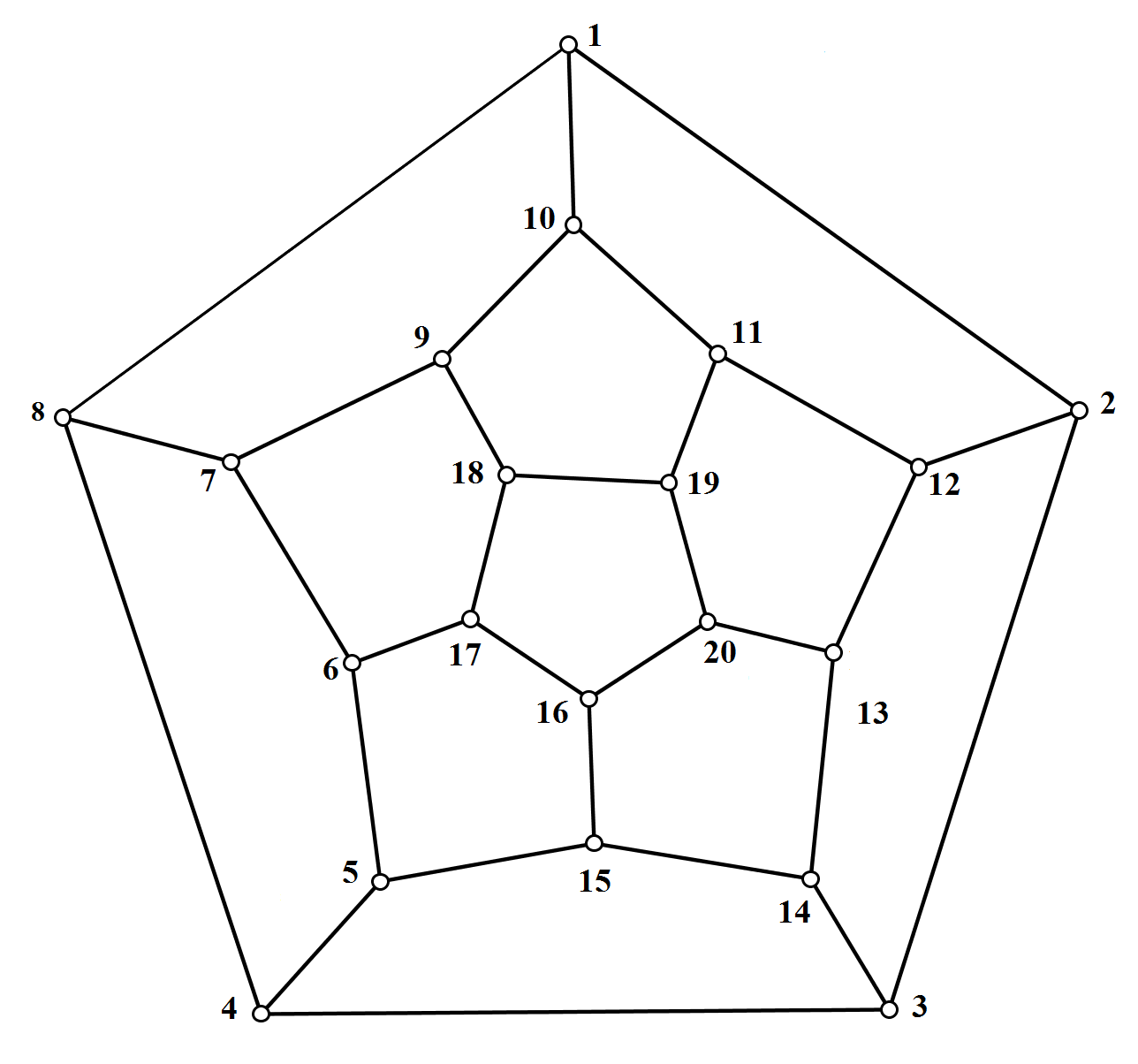}
 \end{center}
\vspace*{-0.1in}\begin{figure}[h]
\centering
{Figure 3. Dodecahedron}
\label{figA.1.1}
\end{figure}
\end{example}

The  next  table contains  C-$\delta$  graphs  $G$  taken  from  \cite{RW}  and upper bounds for   $\msr(\overline{G})$ given  by  $\Delta(G)+1$ are  found in \cite{PD}.\\

\begin{table}[h]
\centering
\caption{Table  of  C-$\delta$  graphs  $G$  taken  from  \cite{RW}  and upper bounds for   $\msr(\overline{G})$ given  by  $\Delta(G)+1$.  }\label{App2}
\end{table}

\begin {center}
\begin{tabular}{|c|c|c|c|c|}
\hline
Family & Name of  Graph      && $\msr(\overline{G})$\\
&$G$&$\vert G\vert$&$\le  \Delta(G)+1 $\\ \hline
Archimedean Graphs&&&\\ \hline
&Cuboctahedron&12&$4$\\ \hline
&Icosidodecahedron&30&$5$\\ \hline
&Rhombicuboctahedron&24&$5$\\ \hline
&Rombicosidodecahedrom&60&$6$\\ \hline
& Snub cube&24&$6$\\ \hline
& Snub dodecahedrom&60&$6$\\ \hline
&Truncated cube&24&$4$\\ \hline
&Truncated Cuboctahedron&48&$4$\\ \hline
&$G$&$\vert G\vert$&$\le  \Delta(G)+1 $\\ \hline
&Truncated dodecahedron&60&$4$\\ \hline
&Truncated icosahedrom &60&$4$\\ \hline
& Truncated icosidodecahedrom&120&$6$\\ \hline
&Truncated Tetrahedron&12&$4$\\ \hline
&Truncated octahedron&24&$4$\\ \hline
Antiprisms &$2n, n\in \N, \ n\ge 3$&$2n, n\ge 3$&$5$\\ \hline
&$4$-antiprism&$8$&$5$\\ \hline
&$5$-antiprism&$10$&$5$\\ \hline
Cages&&&\\ \hline
&Balaban's $(3,10)$ cage&70&$4$\\ \hline
&Foster $(5,5)$ cage&30&$6$\\ \hline
&Harries's $(3,10)$ cage &70&$4$\\ \hline
&Headwood's $(3,6)$ cage&14&$4$\\ \hline
&MacGee's $(3,7)$ cage&24&$4$\\ \hline
&Petersen's $(3,5)$ cage&10&4\\ \hline
&Robertson's $(5,5)$ cage&30&$6$\\ \hline
&Robertson's $(4,5)$ cage&19&$5$\\ \hline
&The Harries-Wong $(3,10)$ cage&70&$4$\\ \hline
&The $(4,6)$ cage&26&$5$\\ \hline
&Tutte's$(3,8)$ cage&30&$4$\\ \hline
&Wongs's $(5,5)$ cage &30&$4$\\ \hline
&The Harries-Wong $(3,10)$ cage&70&$4$\\ \hline
&The $(4,6)$ cage&26&$5$\\ \hline
&Tutte's$(3,8)$ cage&30&$4$\\ \hline
&Wongs's $(5,5)$ cage &30&$4$\\ \hline

\end{tabular}
\end{center}
\newpage
\begin{center}
\begin{tabular}{|c|c|c|c|c|}
\hline
Family & Name of  the   Graph      && $\msr(\overline{G})$\\
&$G$&$\vert G\vert$&$\le  \Delta(G)+1 $\\ \hline
Blanusa Snarks&&&\\ \hline
&Type 1:  26 vertices&26&4\\ \hline
&Type 2: 26 vertices&26&4\\ \hline
&Type 1: 34 vertices&34&4\\ \hline
&Type 2: 34 vertices&34&4\\ \hline
&Type 1: 42 vertices&42&4\\ \hline
&Type 2: 34 vertices&42&4\\ \hline
Generalized Petersen Graphs&&&\\ \hline
&Gp1&10&$4$\\ \hline
&Gp2&12&$4$\\ \hline
&Gp3&14&$4$\\ \hline
&Gp4&16&$4$\\ \hline
&Gp5&16&$4$\\ \hline
&Gp6&18&$4$\\ \hline
&Gp7&18&$4$\\ \hline
&Gp8&20&$4$\\ \hline
&Gp9&20&$4$\\ \hline
&Gp10&20&$4$\\ \hline
&Gp11&22&$4$\\ \hline
&Gp12&22&$4$\\ \hline
&Gp13&24&$4$\\ \hline
&Gp14&24&$4$\\ \hline
&Gp15&24&$4$\\ \hline
&Gp16&24&$4$\\ \hline
Non-Hamiltonian Cubic &&&\\ \hline
&Grinberg's Graph&44&4\\ \hline
&Tutte's Graph&46&$4$\\ \hline
&(38 vertices)&38&4\\ \hline
&(42 vertices)&42&$4$\\ \hline
Platonic Graphs&&&\\ \hline
&Cube&8&$4$\\ \hline
&Dodecahedron&20&$4$\\ \hline
Prisms &$n$-prism,$n\ge 4$ &$2n$&$4$\\ \hline
&$4$-prism &$8$&$4$\\ \hline
&$5$-prism &$10$&$4$\\ \hline
Snarks&&&\\ \hline
&Celmins-Swarf snark 1&26&$4$\\ \hline
&Celmins-Swarf snark 2&26&$4$\\ \hline
&Double Star snark&30&$4$\\ \hline
&Flower snark $J_7$&28&$4$\\ \hline
&Flower snark $J_9$&36&$4$\\ \hline
&Flower snark $J_{11}$&44&$4$\\ \hline

\end{tabular}
\end{center}
\newpage
\begin{center}
\begin{tabular}{|c|c|c|c|c|}
\hline
Family & Graph      && $\msr(\overline{G})$\\
&$G$&$\vert G\vert$&$\le  \Delta(G)+1 $\\ \hline
&Hypercube &$2^4$&$5$\\ \hline
&Loupekine's snark 1 (Sn28)&22&$4$\\ \hline
&Loupekine's snark 2 (Sn29)&22&$4$\\ \hline
&The Biggs-Smith&102&$4$\\ \hline
&The Greenwood-Gleason&16&$6$\\ \hline
&The Szekeres snark&50&$4$\\ \hline
&Watkin's snark &50&$4$\\ \hline
Miscelaneous Regular Graphs&&&\\ \hline
&Chvatal's graph&12&$5$\\ \hline
&Cubic Graph with no perfect matching &$16$&$4$\\ \hline
&Cubic Identity  graphs &$12$&$4$\\ \hline
&Folkman's graph&20&$5$\\ \hline
&Franklin's graph&12&$4$\\ \hline
&Herschel's graph&11&$5$\\ \hline
&Hypercube &$16$&$4$\\ \hline
&Meredith's graph&70&$4$\\ \hline
&Mycielslski's graph&11&$6$\\ \hline
&The  Greenwood-Gleason graph &$16$&$6$\\ \hline
&The Goldner-Harary dual&&$$\\
&( the truncated Prism)&18&$4$\\ \hline
&Tietze's graph&11&$4$\\ \hline
\end{tabular}
\end{center}

\section{Proof  of  Delta Conjecture}
In  this  section  we  give  an argument   which   prove  that  Delta  Conjecture  is  true  for any simple graph  not  necessarily connected  as a  generalization of  the  result  given  in \cite{PD}.   For  that  purpose  We    define  a  generalization  of  C-$\delta$ graphs called  {\bf  extended C-$\delta$ g graph}.

Previously,  we stablish that    Delta  Conjecture  holds for  $\delta-$Graphs. The  condition $2\le\Delta(G)\le\vert G\vert-2$ in the proof  of \ref{main} was given as  a  sufficient  condition  to obtain  that  the  graph  complement of  a C-$\delta$ graph is connected. We  will see  that the  condition  of connectivity  of  a  C-$\delta$  graphs   is not  necessary in order to proof  Delta  Conjecture  when   using  the  result \ref{main}.

Hence,  we can  define a  generalization  of  C-$\delta$  graphs  in   the  following  way.

\defn{ A {\bf extended C-$\delta$ graph} $G'$  is a  simple  graph  which  is  the  disjoint union  of  a  simple connected graph $G, \vert G\vert\ge 4$ ( not  necessarily connected) and a path $P_n$,  where  $n= 2\Delta(G)+2, n\ge 4$.  That  is
$$
G'=P_n\sqcup G ; n=2\Delta(G)+2.
$$}

All  vertices  of $G$ are connected  with  all vertices of $P_n$ in  $\overline{G'}$.  As  a  consequence $\overline{G'}$ is a simple  and  connected graph.

\newpage
\defn {A  graph $G, \vert G\vert\ge 4$ has  a {\bf C-$\delta$  construction} if  it  can  be  constructed   starting    with $K_3$ or  $P_3$  and  by  adding  one  vertex  at a time  in  such  a  way  that  the  newest  vertex $v_m, m\ge 4$ is adjacent to  at most $\dis{\left\lfloor\frac{m}{2}-1\right\rfloor}$ of  the prior  vertices $v_1,v_2,\dots,v_{m-1}$}.

\begin{center}
\includegraphics[height=50mm]{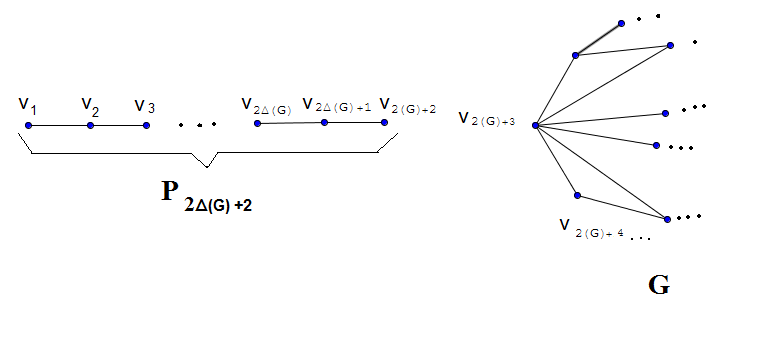}
 \end{center}
\vspace{-0.1in}\begin{figure}[h]
\centering
\caption{Extended C-$\delta$ Graph }\label{fig1}
\end{figure}
\proposition\label{propo1}{ Let  $G'=P_n\sqcup G$ be  an  extended C-$\delta$  graph. Then  $\overline{G'}$  has an  orthogonal  representation  in $\Re^{\Delta(G)+1}$. }

\noindent\MakeUppercase{Proof}:

Let $G'(V',E')$  be an  extended  C-$\delta$ graph.  Then $G'=P_n\sqcup G ; n=2\Delta(G)+2$  and  $G(V,E)$  is a  simple  graph. Since  $P_n$  is a  C-$\delta$ graph we  know   that  we can label  its  vertices  in  such  a  way  that if $v_2,\dots,v_n$  are  its   vertices   then  $\overrightarrow{v_1},\dots,\overrightarrow{v_n}$ is an orthogonal  representation of  its  graph  complement $\overline{P_n}$ in  $\Re^3$. But  since $\Delta(G)\ge 2 $  because   $G$ is  connected   and  $\vert G\vert\ge 4$ then  we  can also obtain  an  orthogonal  representation  of  $P_n$  in $\Re^{\Delta(G)+1}$ getting $2\Delta(G)+2$ vectors    for $\overline{G'}$ in $\Re^{\Delta(G)+1}$ using  the  C-$\delta$  construction.

Thus,  in  $G'$, $v_{2\Delta(G)+2}$ is   adjacent with  all  prior  vertices $v_1,\dots,v_{2\Delta(G)+1}$ but  at most $$\dis{\left\lfloor\frac{2\Delta(G)+2}{2}-1\right\rfloor}\\ =\Delta(G)\ge 2$$
vertices. Actually   to  all of  them but  one.

Now,  choose a  vertex $v'$  in  $\overline{G}$  and label it  as $ v'=v_{2\Delta(G)+3}$. In  $\overline{G'}$ $v_{2\Delta(G)+3}$ is  adjacent
with  all of  the  vertices   of  $P_n$. As a consequence, $v_{2\Delta(G)+3}$ satisfies  the  delta  construction   in $G'$.

If  $Y_{2\Delta(G)+2}$  is  the induced  graph  of $G'$  given  by  $v_1,\dots,v_{2\Delta(G)+2}$  then  $Y_{2\Delta(G)+3}=Y_{2\Delta(G)+2}\cup \{v_{2\Delta(G)+3}\}$ is  simple  and  connected  and  $\overline{Y}_{2\Delta(G)+3}$ can  be  constructed   by using  $\delta$-construction because.
$v_{2\Delta(G)+3}$ is  adjacent with  all previous  vertices  $v_1,v_2,\dots v_{2\Delta(G)+2}$ in $\overline{G'}$. Then it is  adjacent   with all previous   vertices in $Y_{2\Delta(G)+2}$ but  at most $\dis{\left\lfloor \frac{2\Delta(G)+3}{2}-1\right\rfloor}\ge \Delta(G)$.

Now, by  labeling  the  remaining  vertices  in  $G'$ which  are    vertices  in  G in  any  random sequence    to  obtain  $v_{2\Delta(G)+4,\dots,v_{2\Delta(G)+2+\vert G\vert}}$  we  get     a sequence  of  induced  subgraph  of  $\overline{G}$
$$
Y_{2\Delta(G)+4}\subseteq Y_{2\Delta(G)+5}\subseteq\dots\subseteq Y_{2\Delta(G)+2+\vert G\vert}=\overline{G'}
$$

All  of  these   induced  subgraphs  can  be  constructed  using  $\delta$-construction.  As  a consequence $Y_{2\Delta(G)+2+\vert G\vert}\\ =\overline{G'}$  can  be  constructed using $\delta$-construction   which  implies   that    there  is  an  orthogonal  representation  $\overrightarrow{v}_1,\overrightarrow{v}_2,\dots,\overrightarrow{v}_{2\Delta(G)+2+\vert G\vert}$ of  the vertices  of  $\overline{G'}$ at  $\Re^{\Delta(G')+1}$.

But $\Delta(G')\ge \Delta(G)$ since  $\vert G\vert\ge 4$,  $\overline{G}$ is  simple a nd  connected, and  $G$   is an induced  graph  of  $G'$.

Then    we  can  get  the  orthogonal  representation of  $\overline{G'}$ in  $\Re^{\Delta(G)+1}$.

Finally,  if $\overrightarrow{v}_1,\overrightarrow{v}_2,\dots,\overrightarrow{v}_{2\Delta(G)+2+\vert G\vert}$ is   the orthogonal  representation of  $G'$ in  $\\Re^{\Delta(G)+1}$   take  the  vectors $\overrightarrow{v}_{2\Delta(G)+3},\overrightarrow{v}_{2\Delta(G)+4},\dots,\overrightarrow{v}_{2\Delta(G)+2+\vert G\vert}$.  These  vectors  satisfy  all the  adjacency  conditions  and  orthogonal  conditions  of  $\overline{G}$  because $\overline{G}$ is an  induced  subgraph  of  $\overline{G'}$. As  a consequence, $\overrightarrow{v}_{2\Delta(G)+3},\overrightarrow{v}_{2\Delta(G)+4},\dots,\overrightarrow{v}_{2\Delta(G)+2+\vert G\vert}$ is  an orthogonal   representation  of $\overline{G}$ in  $\Re^{\Delta(G)+1}$.\\ .\hfill$\square$

\theorem\label{teo1}{If  $G$ is  a  simple  connected  graph, $\vert G\vert \ge 4$  then  $G$  satisfies  Delta  conjecture}.

\noindent\MakeUppercase{Proof}:

Let $G$  be  a  simple  connected graph. Since  $G$ can  be seen  as a  component   of   a extended  C-$\delta$ graph $G'= P_{2\Delta(G)+2} \sqcup G$  by  the  proposition   proved  above $G$   has  an  orthogonal  representation   in   $\Re^{\Delta(G)+1}=\Re^{\vert G\vert-\delta(G)}$  which   implies   that
$\msr(G)\le \vert G\vert-\delta(G)$ . As  a  consequence   delta  conjecture  holds    for any  simple connected  graph  $G$ with  $\vert G\vert\ge 4$.\hfill $\square$

Finally, by  using   extended C-$\delta$ graphs $G'=P_{2\Delta(G)+2}\sqcup G$   for  all  $G, \vert G\vert \le 3$ and   the  technique   described  in  the proof  of   the proposition  above  or  any other    way  it is  easy    to  check that    all   of  simple  connected  graphs    with  $\vert G\vert \le 3$   satisfies Delta  conjecture.  As  a consequence   we  have the  following  theorem:

\corollary{ Let $G$  be a  simple  and  connected  graph. Then $G$   satisfies delta  conjecture}

\noindent\MakeUppercase{Proof}:

From \ref{teo1} we know  that  delta  conjecture  hold   for  any  simple  graph $G, \vert G\vert \ge 4$.  Checking  all  cases for   all  simple  connected  graph $G, \vert G\vert \le 3 $  we complete  the proof  for  delta conjecture. \hfill $\square$

\section{Conclusion}

In  this paper  we  proved  the  delta  conjecture  as  a  main  result. Also   we applied  the technique for finding  the  minimum  semidefinite  rank of  a  C-$\delta$  to  give a  table  of  upper bounds   of  a  large amount  of  families  of   simple  connected  graphs.  These upper  bounds  will  be usefull  in the  study   of the minimum  semidefinite  rank  of  a  graph.

In  the  future, the  techniques applied in  this paper   could  be  useful  to  solve  other  problems   related  with   simple connected  graphs and   minimum  semidefinite  rank.
\section{Acknowledment}
I  would like  to  thanks   to  my  advisor  Dr. Sivaram Narayan  for   his  guidance and  suggestions  of  this  research.  Also  I  want  to  thank  to  the math  department  of  University  of  Costa  Rica and  Universidad  Nacional Estatal  a Distancia  because  their  sponsorship  during my  dissertation  research  and    specially  thanks  to   the math  department    of  Central Michigan University   where  I  did  the  researh   for  C-$delta$ graphs which  was  a paramount  research  to proof delta  conjecture.
\renewcommand{\baselinestretch}{1} 

\end{document}